\documentclass[a4paper,11pt]{article}
\usepackage{amsmath, amssymb, amsthm}
\usepackage{graphicx} 

\begin{document}

\def\endproof{\hfill$\square$}

\def\N{{\mathbf N}} \def\f{\varphi} \def\Im{{\rm Im}}
\def\lg{{\rm lg}} \def\Sym{{\rm Sym}} \def\Id{{\rm Id}}
\def\mod{{\rm mod}} \def\Ker{{\rm Ker}}\def\Conj{{\rm Conj}}
\def\Z{{\mathbb Z}}\def\sgn{{\rm sgn}} \def\Aut{{\rm Aut}}
\def\AA{{\mathcal A}} \def\BB{{\mathcal B}} \def\PP{{\mathcal P}}
\def\R{{\mathbb R}} \def\C{{\mathbb C}}

\title{\bf{On a theorem of Artin, II}}
 
\author{\textsc{Nuno Franco and Luis Paris}}

\date{\today}

\maketitle

\begin{abstract}
\noindent 
This paper is a sequel of \cite{CoPa}. Let $A$ be an Artin group, let $W$ be its associated 
Coxeter group, and let $CA$ be its associated coloured Artin group, that is, the kernel of the 
standard epimorphism $\mu: A \to W$. We determine the homomorphisms $\f: A \to W$ that verify 
$\Im \f \cdot Z(W)= W$, for $A$ irreducible and of spherical type, and we prove that $CA$ is a 
characteristic subgroup of $A$, if $A$ is of spherical type but not necessarily irreducible.
\end{abstract}

\noindent
{\bf AMS Subject Classification:} Primary 20F36. Secondary 20F55.

\section{Introduction}

A {\it Coxeter matrix} of rank $n$ is a square matrix $M=(m_{i\,j})_{1 \le i,j\le n}$ such 
that $m_{i\,i}=1$ for all $i=1, \dots, n$, and $m_{i\,j}=m_{j\,i}\in \{2,3,4, \dots, +\infty\}$ 
for all $i,j \in \{1, \dots, n\}$, $i \neq j$. A Coxeter matrix $M=(m_{i\,j})$ is usually 
represented by its {\it Coxeter graph}, $\Gamma$, which is defined as follows. The set of 
vertices of $\Gamma$ is $\{1, \dots, n\}$, two vertices $i,j$ are joined by an edge if $m_{i\,j} 
\ge 3$, and this edge is labeled by $m_{i\,j}$ if $m_{i\,j}\ge 4$.

\bigskip\noindent
Let $\Gamma$ be a Coxeter graph with set of vertices $\{1, \dots, n\}$. For two objects $a,b$ and 
$m \in \N$, we define the word
\[
\omega(a,b:m)= \left\{
\begin{array}{ll}
(ab)^{{m \over 2}}&\quad\text{if }m\text{ is even}\,,\\
(ab)^{{m-1\over 2}}a&\quad\text{if }m\text{ is odd}\,.\\
\end{array}\right.
\]
The {\it Artin group of type $\Gamma$} is the group $A=A(\Gamma)$ generated by $\sigma_1, \dots, 
\sigma_n$, and subject to the relations
\[
\omega( \sigma_i, \sigma_j: m_{i\,j}) = \omega( \sigma_j, \sigma_i: m_{i\,j}) \quad \text{for }
1 \le i<j\le n \text{ and } m_{i\,j}< +\infty\,,
\]
where $M=(m_{i\,j})$ is the Coxeter matrix of $\Gamma$. The {\it Coxeter group of type $\Gamma$} 
is the group $W=W(\Gamma)$ generated by $s_1, \dots, s_n$, and subject to the relations
\begin{gather*}
s_i^2=1 \quad \text{for } 1 \le i\le n\,,\\
(s_is_j)^{m_{i\,j}} =1 \quad \text{for } 1\le i<j\le n \text{ and }m_{i\,j}<+\infty\,.
\end{gather*}
Note that $W$ is the quotient of $A$ by the relations $\sigma_i^2=1$, $1 \le i\le n$, and $s_i$ 
is the image of $\sigma_i$ under the quotient epimorphism. Throughout the paper, this epimorphism 
will be denoted by $\mu: A \to W$, and called {\it standard epimorphism}.

\bigskip\noindent
The number $n$ of generators is called the {\it rank} of the Artin group (and of the Coxeter 
group). We say that $A$ is {\it irreducible} if $\Gamma$ is connected, and that $A$ is of {\it 
spherical type} if $W$ is finite. Define the {\it coloured Artin group} $CA=CA(\Gamma)$ to be the 
kernel of the standard epimorphism $\mu: A \to W$.

\bigskip\noindent
Artin groups were first introduced by Tits \cite{Tit2} as extensions of Coxeter groups. Later, 
Brieskorn \cite{Bri} gave a topological interpretation of the Artin groups and
coloured Artin groups of spherical type in terms of 
regular orbit spaces as follows.

\bigskip\noindent
Define a {\it finite reflection group} to be a finite subgroup $W$ of $O(n,\R)$ generated by 
reflections, where $n$ is some positive integer. A classical result due to Coxeter \cite{Cox1}, 
\cite{Cox2}, states that $W$ is a finite reflection group if and only if $W$ is a finite Coxeter 
group. Assume this is the case. Let $\AA(\Gamma)$ be the set of reflecting hyperplanes of 
$W=W(\Gamma)$, and, for $H \in \AA(\Gamma)$, let $H_\C$ denote the hyperplane of $\C^n$ having 
the same equation as $H$. Let
\[
M(\Gamma)= \C^n \setminus \left( \bigcup_{H \in \AA( \Gamma)} H_\C \right) \,.
\]
Then $M(\Gamma)$ is a connected submanifold of $\C^n$, the group $W(\Gamma)$ acts freely on 
$M(\Gamma)$, and the quotient $N(\Gamma)= M(\Gamma)/ W(\Gamma)$ is isomorphic to the complement 
in $\C^n$ of an algebraic variety called {\it discriminental variety} of type $\Gamma$. By 
\cite{Bri}, the fundamental group of $N(\Gamma)$ is the Artin group $A(\Gamma)$, 
the fundamental group of $M(\Gamma)$ is the associated coloured Artin group $CA(\Gamma)$, and the 
exact sequence associated to the regular covering $M(\Gamma) \to N(\Gamma)$ is $1 \to CA(\Gamma) 
\to A(\Gamma) \xrightarrow{\mu} W(\Gamma) \to 1$.

\bigskip\noindent
The set $\AA(\Gamma)$ is called {\it Coxeter arrangement}. Coxeter arrangements have been widely 
studied from several points of view (combinatorics, topology, homology, ...). For instance, by 
\cite{Del}, the space $M(\Gamma)$ is a classifying space for $CA(\Gamma)$. The cohomology of 
$M(\Gamma)$ (and, hence, of $CA(\Gamma)$) has been calculated in \cite{Bri2}. In particular, the 
Poincar\'e polynomial of $CA(\Gamma)$ is $\prod_{j=1}^n (1+m_jt)$, where $m_1, \dots, m_n$ are the 
exponents of $W(\Gamma)$. We refer to \cite{OrTe} for a detailed exposition on arrangements of 
hyperplanes which includes a good account on the theory of Coxeter arrangements.

\bigskip\noindent
In the case $\Gamma=A_{n-1}$, the group $W(\Gamma)$ is the symmetric group $\Sym_n$, $A(\Gamma)$ 
is the braid group $\BB_n$ introduced by Artin in 1925 \cite{Art1}, $CA(\Gamma)$ is the pure braid 
group $\PP\BB_n$, and $N(\Gamma)$ is the space of configurations of $n$ (unordered) points in $\C^n$.

\bigskip\noindent
Call two homomorphisms $\f_1, \f_2: A \to W$ {\it equivalent} if there is an automorphism 
$\alpha$ of $W$ such that $\f_2= \alpha \circ \f_1$.

\bigskip\noindent
In 1947, Artin \cite{Art2} determined all the epimorphisms $\BB_n \to \Sym_n$, up to equivalence, 
and proved that the pure braid group $\PP\BB_n$ is a characteristic subgroup of $\BB_n$. This 
last result is of importance, for example, in the calculation of the automorphism group of 
$\BB_n$ (see \cite{DyGr}). The classification of the epimorphisms $A \to W$, for $A$ of type 
$B_n$ and $D_n$, was almost done by Lin in \cite{Lin}. The complete classification of the 
epimorphisms $A \to W$, for $A$ of spherical type and irreducible, was done later in \cite{CoPa}. 
It is also proved in \cite{CoPa} that $CA$ is a characteristic subgroup of $A$, in the case when 
$A$ is of spherical type and irreducible.

\bigskip\noindent
The goal of the present paper is to extend this last result to the non-irreducible spherical type 
Artin groups, that is, we prove that $CA$ is a characteristic subgroup of $A$, when $A$ is
of spherical type, but not necessarily irreducible. 
This extension is not immediate as we need extra arguments to palliate the case by 
case treatment made in \cite{CoPa}. Moreover, one can believe that the automorphism group of a 
non-irreducible Artin group should be much larger that the automorphism group of an irreducible 
one. For instance, $\Z$ is an irreducible Artin group, while $\Z^n$ is not irreducible if $n \ge 
2$. The outer automorphism group of the braid group is of order $2$ (see \cite{DyGr}), but it is 
not true that the outer automorphism group of an irreducible (spherical type) Artin group is 
always finite (see \cite{ChCr}).

\bigskip\noindent
In order to achieve our main result, we first need to classify the epimorphisms $A \to W$ ``up 
to the center'', for all irreducible and spherical type Artin groups. This is the object of 
Section 2. In Section 3 we prove that, for $A$ irreducible, if $\Phi: A \to A$ is a homomorphism 
such that $\mu \circ \Phi: A \to W$ is an epimorphism ``up to the center'', then $\Phi$ must 
send $CA$ into $CA$. Finally, we prove our main result in Section 4.

\bigskip\noindent
From now on, $\Gamma$ denotes a spherical type Coxeter graph, $A=A(\Gamma)$ denotes the Artin 
group associated to $\Gamma$, and $W=W(\Gamma)$ denotes the Coxeter group associated to 
$\Gamma$.

\bigskip\noindent
{\bf Acknowledgments.} We would like to thank the {\it Centre de Calcul MEDICIS}, at the research
center Sciences et Technologies de l'Information et de la Communication at X - FRE CNRS 2341,
the {\it Centopeia}, at the Centro de Fsica Computacional da Universidade de Coimbra, and to
{\it CIMA-UE}, Universidade de \'Evora, for allowing the access to their computers.

\section{Up-to-center-epimorphisms}

Let $Z(W)$ denote the center of $W$. Call a homomorphism $\f: A \to W$ an {\it up-to-center-epimorphism} 
if $W= \Im\f \cdot Z(W)$. Note that, if $W$ is centerless, then any up-to-center-epimorphism 
is an epimorphism. Conversely, any epimorphism is an up-to-center-epimorphism. An 
up-to-center-epimorphism is said to be {\it proper} if it is not an epimorphism.

\bigskip\noindent
An up-to-center-epimorphism $\f: A \to W$ is called {\it ordinary} if $\f(\sigma_i^2)=1$ for all 
$i=1, \dots, n$. The up-to-center-epimorphisms that are not ordinary are called {\it 
extraordinary}. Note that an up-to-center-epimorphism is ordinary if and only if $CA$ lies in its 
kernel. Note also that an epimorphism is ordinary if and only if it is equivalent to the standard 
one.

\bigskip\noindent
The aim of the present section is the classification of the 
up-to-center-epimor\-phisms $A \to W$, up to equivalence, for all irreducible and spherical type 
Artin groups. 
Note that the classes of epimorphisms have been already classified in \cite{CoPa}, so our work 
essentially consists on determining the classes of proper up-to-center-epimorphisms.

\bigskip\noindent
We start with some preliminary well-known results on Coxeter groups and Artin groups.

\bigskip\noindent
The connected spherical type Coxeter graphs are precisely the Coxeter graphs $A_n$ ($n \ge 1$), 
$B_n$ ($n\ge 2$), $D_n$ ($n \ge 4$), $E_6$, $E_7$, $E_8$, $F_4$, $H_3$, $H_4$, $I_2(p)$ ($p\ge 
5$) illustrated in Figure 1. Here we use the notation $I_2(6)$ for the Coxeter graph $G_2$. We 
may also use the notation $I_2(3)$ for $A_2$, and $I_2(4)$ for $B_2$. We label the vertices of 
each of these graphs as shown in Figure 1.

\begin{figure}
\begin{center}
\includegraphics[width=8cm]{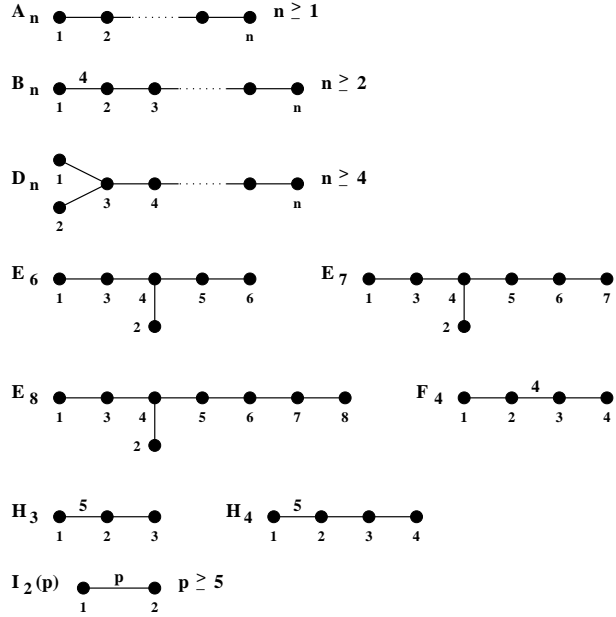}
\end{center}
\caption{The connected spherical type Coxeter graphs.}
\end{figure}

\bigskip\noindent
For $w \in W$, we denote by $\lg(w)$ the word length of $w$ with respect to $S$. The group $W$ 
has a unique element of maximal length, $w_0$, which satisfies $w_0^2=1$ and $w_0 S w_0=S$. The 
standard epimorphism $\mu: A \to W$ has a natural set-section $T:W \to A$ defined as follows. Let 
$w \in W$, and let $w=s_{i_1} s_{i_2} \dots s_{i_l}$ be a reduced expression for $w$ ({\it i.e.} 
$l=\lg (w)$). Then $T(w)= \sigma_{i_1} \sigma_{i_2} \dots \sigma_{i_l} \in A$. By Tits' solution 
to the word problem for Coxeter groups (see \cite{Tit}), the definition of $T(w)$ does not depend 
on the choice of the reduced expression. Define the {\it fundamental element} of $A$ to be 
$\Delta= T(w_0)$, where $w_0$ is the longest element of $W$.

\bigskip\noindent
For a group $G$ and a subset $X \subset G$, we denote by $Z(G)$ the center of $G$, and by 
$N_G(X)$ the {\it normalizer} of $X$ in $G$, that is, $N_G(X)= \{ g\in G; 
gXg^{-1} = X \}$. The following proposition is a mixture of several well-known facts on 
finite Coxeter groups and spherical type Artin groups. The results concerning the Coxeter groups 
can be found in \cite{Bou}, and the ones concerning the Artin groups can be found in 
\cite{BrSa} and \cite{Del}.

\bigskip\noindent
{\bf Proposition 2.1.} {\it Assume $\Gamma$ to be connected. Let $w_0$ be the longest element of 
$W$, and let $\Delta$ be the fundamental element of $A$.
\begin{enumerate}
\item There exists a permutation $\xi \in \Sym_n$ of order 2 such that $w_0 s_i w_0 = s_{\xi(i)}$ 
and $\Delta \sigma_i \Delta^{-1} = \sigma_{\xi(i)}$ for all $i=1, \dots, n$.
\item The group $N_W(\{s_1, \dots, s_n\})$ is cyclic of order 2 generated by $w_0$, and $N_A( \{ 
\sigma_1, \dots, \sigma_n \})$ is an infinite cyclic subgroup generated by $\Delta$.
\item If $\xi=\Id$, then $Z(W)=N_W(\{ s_1, \dots, s_n \})$ is the cyclic subgroup of $W$ of order 
2 generated by $w_0$, and $Z(A)=N_A(\{ \sigma_1, \dots, \sigma_n \})$ is the infinite cyclic 
subgroup of $A$ generated by $\delta= \Delta$. If $\xi \neq \Id$, then $Z(W)= \{1\}$, and $Z(A)$ 
is the infinite cyclic subgroup of $A$ generated by $\delta= \Delta^2$.
\item We have $\xi \neq \Id$ if and only if $\Gamma \in \{ A_n; n \ge 2\} \cup \{D_n; n\ge 5 
\text{ and } n \text{ odd}\} \cup \{I_2(p); p\ge 5 \text{ and }p \text{ odd}\} \cup \{E_6\}$.
\end{enumerate}}

\noindent
The generator $\delta$ of $Z(A)$ in the above proposition will be called the {\it standard 
generator} of $Z(A)$.

\bigskip\noindent
Now, using Proposition 2.1, we can prove the following.

\bigskip\noindent
{\bf Lemma 2.2.} {\it Assume $\Gamma$ to be connected.
\begin{enumerate}
\item There exists a proper up-to-center-epimorphism $A \to W$ if and only if $\Gamma \in \{ 
I_2(p); p\ge 6 \text{ and }p\equiv 2(\mod\, 4)\} \cup \{B_n; n\ge 3 \text{ and }n \text{ odd}\} 
\cup \{A_1, H_3, E_7\}$.
\item If $\f:A \to W$ is a proper up-to-center-epimorphism, then $W= \Im\, \f \times Z(W)$.
\end{enumerate}}

\noindent
{\bf Proof.} We suppose that there exists a proper up-to-center-epimorphism $\f: A \to W$, and 
we write $H=\Im\, \f$. So, $W=H \cdot Z(W)$ and $H \neq W$. Clearly, we must have $Z(W) \neq 
\{1\}$, thus, by Proposition 2.1, $\Gamma \in \{ I_2(p); p\ge 4 \text{ and }p\text{ even}\} \cup 
\{B_n; n\ge 3\} \cup \{ D_n; n\ge 4 \text{ and } n \text{ even}\} \cup \{A_1, E_7, E_8, F_4, H_3, 
H_4 \}$. Moreover, again by Proposition 2.1, we have $|Z(W)|=2$, thus $|W:H| \le 2$. We 
cannot have $|W:H|=1$ since $H \neq W$, therefore $|W:H|=2$ and $W=H \times Z(W)$. 

\bigskip\noindent
The quotient of $W$ by $H$ determines an epimorphism $\gamma: W \to C_2=\{\pm 1\}$ which 
satisfies $\gamma(w_0)=-1$ and $\Ker \gamma = H$. Conversely, if $\gamma: W \to \{\pm 1\}$ is an 
epimorphism which satisfies $\gamma(w_0)=-1$, then $W=\Ker \gamma \times Z(W)$ and $\pi \circ 
\mu: A \to W= \Ker \gamma \times Z(W)$ is a proper up-to-center-epimorphism, where $\pi: W \to 
\Ker \gamma \subset W$ denotes the projection on the first coordinate. So, a proper up-to-center-epimorphism 
exists if and only if there exists an epimorphism $\gamma: W \to C_2$ which 
satisfies $\gamma(w_0)=-1$.

\bigskip\noindent
Recall from \cite{Bou} that $w_0= (s_1s_2 \dots s_n)^{h/2}$, where $h$ is the Coxeter number of 
$W$. Thus, if $h/2$ is even, then $\gamma(w_0)=1$ for any epimorphism $\gamma: W \to C_2$. So, we 
can assume that $h/2$ is odd, that is, $\Gamma \in \{I_2(p); p\equiv 2(\mod\,4)\} \cup \{B_n; 
n\equiv 1(\mod\,2)\} \cup \{A_1,E_7,E_8,H_3,H_4\}$.

\bigskip\noindent
Suppose $\Gamma=I_2(p)$ with $p \equiv 2(\mod\,4)$. Let $\gamma:W \to C_2$ be the epimorphism 
defined by $\gamma(s_1)=-1$ and $\gamma(s_2)=1$. Then $\gamma(w_0)=-1$.

\bigskip\noindent
Suppose $\Gamma=B_n$ with $n \equiv 1(\mod\,2)$. Let $\gamma:W \to C_2$ be the epimorphism 
defined by $\gamma(s_1)=\gamma(s_2) = \dots = \gamma(s_n) = -1$. Then $\gamma(w_0)=-1$.

\bigskip\noindent
Suppose $\Gamma=A_1$. Let $\gamma:W \to C_2$ be the epimorphism defined by $\gamma(s_1)=-1$. Then 
$\gamma(w_0)= \gamma(s_1)= -1$. (This case is obvious.)

\bigskip\noindent
Suppose $\Gamma=H_3$. Let $\gamma:W \to C_2$ be the epimorphism defined by $\gamma(s_1)= 
\gamma(s_2) = \gamma(s_3) =-1$. Then $\gamma(w_0)=-1$.

\bigskip\noindent
Suppose $\Gamma=H_4$. Let $\gamma:W \to C_2$ be an epimorphism. The generator $s_i$ is conjugate 
to $s_1$ for all $i=1, \dots, 4$, thus $\gamma(s_1)= \gamma(s_2)= \gamma(s_3)= \gamma(s_4)$, 
therefore $\gamma(s_1s_2s_3s_4)=1$, hence $\gamma(w_0)=1$.

\bigskip\noindent
Suppose $\Gamma=E_7$. Let $\gamma:W \to C_2$ be the epimorphism defined by $\gamma(s_1)= \dots = 
\gamma(s_7)=-1$. Then $\gamma(w_0)=-1$.

\bigskip\noindent
Suppose $\Gamma=E_8$. Let $\gamma:W \to C_2$ be an epimorphism. Here again, the generator $s_i$ 
is conjugate to $s_1$ for all $i=1, \dots, 8$, thus $\gamma(s_1)= \gamma(s_2)= \dots = 
\gamma(s_8)$, therefore $\gamma(s_1s_2 \dots s_8)=1$, hence $\gamma(w_0)=1$.
\endproof

\bigskip\noindent
As pointed out before, the equivalence classes of epimorphisms $A \to W$ are known (see 
\cite{CoPa}), thus, by the above lemma, we only need to study those Coxeter graphs which lie in 
the set $\{I_2(p); p\equiv 2(\mod\,4)\} \cup \{B_n; n\ge 3\text{ and }n\text{ odd}\} \cup 
\{H_3,E_7\}$. (The case $\Gamma=A_1$ is trivial.) We start with the rank 2 groups.

\bigskip\noindent
{\bf Proposition 2.3.} {\it Suppose $\Gamma=I_2(p)$, where $p \equiv 2(\mod\,4)$. Then there is 
precisely one equivalence class of proper up-to-center-epimorphisms. This class is represented by 
$\mu': A \to W$, where
\[
\mu'(\sigma_1)=s_1\,,\text{ and } \mu'(\sigma_2)= s_2s_1s_2\,.
\]
In particular, since $\mu'$ is ordinary, any proper up-to-center-epimorphism is ordinary.}

\bigskip\noindent
{\bf Proof.} Let $\f: A \to W$ be a proper up-to-center-epimorphism. Write $p=4k+2$. Let 
$\gamma_1, \gamma_2: W \to C_2=\{\pm 1\}$ be the epimorphisms determined by
\[
\begin{array}{cc}
\gamma_1(s_1)=1\,,\quad &\gamma_1(s_2)=-1\,,\\
\gamma_2(s_1)=-1\,,\quad &\gamma_2(s_2)=1\,,
\end{array}
\]
and let $H_1=\Ker \gamma_1$ and $H_2=\Ker \gamma_2$. One can easily verify that $H_1$ and $H_2$ 
are the only subgroups $H$ of $W$ that satisfy $H \times Z(W)=W$, thus either $\Im \f=H_1$ or 
$\Im\f =H_2$. Up to an automorphism of $W$, we can assume that $\Im\f= H_1$.

\bigskip\noindent
Write $t_1=s_1$ and $t_2=s_2s_1s_2$. Then $H_1= \langle t_1,t_2 | t_1^2= t_2^2= (t_1t_2)^{2k+1} 
=1 \rangle \simeq W(I_2(2k+1))$. Consider the epimorphism $\bar\gamma: H_1 \to C_2$ which sends $t_1$ 
and $t_2$ to $-1$.

\bigskip\noindent
We cannot have $\bar\gamma( \f(\sigma_1))= \bar\gamma (\f( \sigma_2))=1$ because, otherwise, we would 
have $\Im\f \subset \Ker(\bar\gamma) \subsetneqq H_1$.

\bigskip\noindent
Suppose that $\bar\gamma( \f( \sigma_1))=-1$ and $\bar\gamma( \f (\sigma_2))=1$. Up to an automorphism of 
$H_1$ (and, hence, of $W=H_1 \times C_2$), we can assume that $\f(\sigma_1)=t_1$ and 
$\f(\sigma_2)= (t_1t_2)^l$ for some $l \in \{0,1, \dots, 2k\}$. Observe that
\begin{multline*}
t_2(t_1t_2)^{l-1} = (t_2(t_1t_2)^{l-1})^{2k+1} = \f( (\sigma_1 \sigma_2)^{2k+1}) \\
= \f( (\sigma_2 \sigma_1)^{2k+1}) = ((t_1t_2)^lt_1)^{2k+1} = (t_1t_2)^lt_1\,,
\end{multline*}
thus $(t_1t_2)^{2l}=1$. This implies that $2k+1$ divides $2l$, hence $l=0$ and $\f (\sigma_2)=1$. 
It follows that $\Im \f = \{t_1,1\} \neq H_1$: a contradiction.

\bigskip\noindent
Similarly, we cannot have simultaneously $\bar\gamma (\f (\sigma_1)) =1$ and $\bar\gamma (\f (\sigma_2))= 
-1$, thus $\bar\gamma( \f (\sigma_1))$ $= \bar\gamma (\f (\sigma_2)) =-1$.

\bigskip\noindent
Let $T= \{wt_1w^{-1}; w\in H_1\} = \{w t_2w^{-1}; w\in H_1\}$. We have $\f (\sigma_1), \f 
(\sigma_2) \in T$, and $\f(\sigma_1), \f (\sigma_2)$ generate $H_1$, thus, up to an automorphism 
of $H_1$, we can suppose that $\f (\sigma_1)=t_1=s_1$ and $\f (\sigma_2)=t_2 =s_2s_1s_2$.
\endproof

\bigskip\noindent
Now, we assume that $\Gamma=B_n$ and $n$ is odd. We have an isomorphism $W\simeq (C_2)^n \rtimes 
\Sym_n$ given by $s_1= [-1,1, \dots, 1] \in (C_2)^n$ and $s_i=(i-1,i) \in \Sym_n$ for all $i=2, 
\dots, n$. The element of maximal length in $W$ is $w_0=[-1,-1, \dots, -1]$ and is central. 
Put $e_i= s_i \dots s_2s_1s_2 \dots s_i = [1, \dots, 1,-1,1, \dots, 1]$ for $i=1, \dots, n$. Then 
$e_1,e_2, \dots, e_n$ are the standard generators of $(C_2)^n$.

\bigskip\noindent
{\bf Proposition 2.4.} {\it Suppose $\Gamma=B_n$ with $n$ odd. There are precisely three classes 
of proper up-to-center-epimorphisms $A \to W$. These are represented by $\mu', \nu_n^3, \nu_n^4$, 
where
\begin{gather*}
\mu'(\sigma_1) =w_0s_1\,,\quad
\mu'(\sigma_i)=s_i 
\text{ for }2\le i\le n\,,\\
\nu_n^3 (\sigma_1) =w_0s_1\,,\quad
\nu_n^3 (\sigma_i)= e_{i-1} s_i
\text{ for }2\le i\le n\,,\\
\nu_n^4 (\sigma_1) =1\,,\quad
\nu_n^4 (\sigma_i) =e_{i-1} s_i
\text{ for } 2\le i\le n\,.
\end{gather*}}

\noindent
{\bf Remark.} The up-to-center-epimorphism $\mu'$ is ordinary, while $\nu_n^3$ and 
$\nu_n^4$ are not.

\bigskip\noindent
{\bf Proof.} Our proof follows the same strategy as the proof of Lemma 3.4 in 
\cite{CoPa}. We consider a proper up-to-center-epimorphism $\f: A \to W$.

\bigskip\noindent
{\bf Assertion 1.} {\it Up to an automorphism of $W$, we can assume that there exist $v_1, v_2, 
\dots, v_n \in (C_2)^n$ such that
$\f (\sigma_1) =v_1$, and $\f(\sigma_i) =v_i \cdot (i-1,i)$ for $2 \le i\le 
n$.}

\bigskip\noindent
{\bf Proof.} For a group $G$ and an element $g \in G$, we denote by $\Conj_g: G \to G$ the inner 
automorphism which sends $h$ to $ghg^{-1}$ for all $h \in G$. Call two homomorphism $\f_1, \f_2: G 
\to H$ {\it conjugate} if there is an $h \in H$ such that $\f_2= \Conj_h \circ \f_2$.

\bigskip\noindent
Let $\eta: A \to \Sym_n$ be the standard epimorphism determined by
$\eta (\sigma_1) =1$, and $\eta (\sigma_i) = (i-1,i)$ for $2 \le i\le n$.
We know by \cite{Zin} that any epimorphism $A \to \Sym_n$ is conjugate to $\eta$. Let $\pi: 
(C_2)^n \rtimes \Sym_n \to \Sym_n$ be the projection on the second coordinate. By the above 
observation, there exists $u \in \Sym_n$ such that $\Conj_u \circ \pi \circ \f = \eta$. Choose 
some $w \in W$ such that $\pi(w)=u$. Then $\eta= \Conj_u \circ \pi \circ \f = \pi \circ \Conj_w 
\circ \f$. This implies that there exist $v_1, v_2, \dots, v_n \in (C_2)^n$ such that $(\Conj_w 
\circ \f)( \sigma_1) = v_1$, and $(\Conj_w \circ \f) (\sigma_i) = v_i \cdot (i-1,i)$ for $2 \le 
i\le n$.

\bigskip\noindent
{\bf Assertion 2.} {\it Up to an automorphism of $W$, we have either
\begin{enumerate}
\item $\f (\sigma_i) = s_i$ for all $2 \le i\le n$; or
\item $\f (\sigma_i) = e_{i-1} s_i$ for all $2\le i\le n$.
\end{enumerate}}

\noindent
{\bf Proof.} Write $v_i= e_1^{a_{i\,1}} e_2^{a_{i\,2}} \dots e_n^{a_{i\,n}} \in (C_2)^n$, where 
$a_{i\,j} \in \Z$, for all $1 \le i\le n$. It is easily checked that the equalities $\f (\sigma_j) 
\f(\sigma_2)= \f (\sigma_2) \f (\sigma_j)$, $4 \le j\le n$, imply
\[
a_{2\,j} \equiv a_{2\,3}\ (\mod\,2) \quad \text{for } 4 \le j\le n\,,
\]
and that the equalities $\f (\sigma_i) \f(\sigma_{i+1}) \f (\sigma_i) = \f (\sigma_{i+1}) \f 
(\sigma_i) \f (\sigma_{i+1})$, $2\le i\le n-1$, imply
\begin{gather*}
a_{i\,j} \equiv a_{2\,3}\ (\mod\,2) \quad \text{for } 2\le i\le n \text{ and } j \neq i,i-1\,,\\
a_{i\,i-1} + a_{i\,i} \equiv a_{2\,1} + a_{2\,2}\ (\mod\,2) \quad \text{for } 2\le i\le n\,.
\end{gather*}
Put $a=a_{2\,3}$. Then $v_i=e_1^a \dots e_{i-2}^a e_{i-1}^{a_{i\,i-1}} e_i^{a_{i\,i}} e_{i+1}^a 
\dots e_n^a$.

\bigskip\noindent
Suppose $a \equiv 1(\mod\,2)$. Consider the automorphism $\alpha_0: W \to W$ defined by
\[
\alpha_0 (s_1)=s_1\,, \quad \alpha_0(s_i)= w_0s_i \text{ for } 2 \le i\le n\,.
\]
Observe that $\alpha_0(u)=u$ for all $u \in (C_2)^n$, thus
\[
(\alpha_0 \circ \f)(\sigma_i) = (e_{i-1}^{a_{i\,i-1}+1} e_i^{a_{i\,i}+1}) \cdot s_i \quad 
\text{for }2\le i\le n\,.
\]
So, we can suppose $a=0$, that is, $\f (\sigma_i) = e_{i-1}^{a_{i\,i-1}} e_i^{a_{i\,i}} s_i$ for 
all $2 \le i\le n$.

\bigskip\noindent
Suppose $a_{n\,n}= \dots = a_{i+1\,i+1} =0$ and $a_{i\,i}=1$. Observe that
\[
(\Conj_{e_{i-1}} \circ \f) (\sigma_j) = \left\{
\begin{array}{ll}
e_{j-1}^{a_{j\,j-1}} e_j^{a_{j\,j}} s_j \quad&\text{if } j \neq i-1,i\,,\\
e_{j-1}^{a_{j\,j-1}+1} e_j^{a_{j\,j}+1} s_j \quad&\text{if }j=i-1 \text{ or }i\,.
\end{array}\right.
\]
So, we can assume that $a_{n\,n} = \dots = a_{i+1\,i+1} = a_{i\,i} =0$. An iteration of this 
argument shows that we can assume that $a_{i\,i}=0$ for all $2\le i\le n$.

\bigskip\noindent
Assume $a_{2\,1}= a_{2\,1} + a_{2\,2} \equiv 0(\mod\,2)$. As $a_{i\,i-1} = a_{i\,i-1} + a_{i\,i} 
\equiv a_{2\,1} + a_{2\,2} (\mod\,2)$, it follows that $a_{i\,i-1} \equiv 0 (\mod\,2)$ and 
$\f(\sigma_i) =s_i$ for all $2\le i\le n$. Suppose $a_{2\,1}= a_{2\,1} +a_{2\,2} \equiv 
1(\mod\,2)$. Then $a_{i\,i-1} \equiv 1(\mod\,2)$ and $\f (\sigma_i)= e_{i-1} s_i$ for all $2\le 
i\le n$.

\bigskip\noindent
{\bf Assertion 3.} {\it There exist $a,b \in \{0,1\}$ such that $\f (\sigma_1)= w_0^a s_1^b$.}

\bigskip\noindent
{\bf Proof.} Write $\f (\sigma_1)= v_1= e_1^{a_1} e_2^{a_2} \dots e_n^{a_n}$, where $a_i \in \Z$. 
Observe that the equalities $\f (\sigma_i) \f (\sigma_1)= \f(\sigma_1) \f(\sigma_i)$, $3 \le i\le 
n$, imply that $a_i \equiv a_2(\mod\,2)$ for all $2 \le i\le n$. Put $a=a_2$ and $b=a_1+a_2$. 
Then $\f (\sigma_1)= w_0^a e_1^b= w_0^a s_1^b$.

\bigskip\noindent
{\bf End of the proof of Proposition 2.4.} First, assume that $\f(\sigma_i) = s_i$ for all $2\le 
i\le n$ (see Assertion 2). If $\f (\sigma_1)=1$, then $\Im \f = \Sym_n$, which is of index $2^n$ 
in $W$, thus $\f$ cannot be an up-to-center-epimorphism. If $\f(\sigma_1)=w_0$, then $\Im\f = 
Z(W) \times \Sym_n$ which is of index $2^{n-1}$ in $W$, thus $\f$ cannot be an up-to-center-epimorphism, 
neither. If $\f(\sigma_1)=s_1$, then $\f= \mu$ (which is not proper), and, if 
$\f(\sigma_1)= w_0s_1$, then $\f= \mu'$.

\bigskip\noindent
Now, assume that $\f(\sigma_i)= e_{i-1} s_i$ for all $2 \le i\le n$. If $\f(\sigma_1)=s_1$, then 
$\f$ is the epimorphism $\nu_n^1$ of \cite{CoPa}, Lemma 3.4, and, if $\f(\sigma_1)= w_0$, then 
$\f$ is the epimorphism $\nu_n^2$ of \cite{CoPa}, Lemma 3.4. If $\f (\sigma_1)= w_0s_1$, then $\f= 
\nu_n^3$, and, if $\f(\sigma_1)= 1$, then $\f= \nu_n^4$.
\endproof

\bigskip\noindent
It remains to investigate the Artin groups of type $H_3$ and $E_7$. 
Recall that two homomorphisms $\f_1, \f_2: A \to W$ are called {\it conjugate} if there exists $w 
\in W$ such that $\Conj_w \circ \f_1=\f_2$. We start with the classification of the conjugacy 
classes of proper up-to-center-epimorphisms $A(H_3) \to W(H_3)$.

\bigskip\noindent
{\bf Lemma 2.5.} {\it Assume $\Gamma=H_3$. Then there are precisely four conjugacy classes of 
proper up-to-center-epimorphisms $A \to W$. These are represented by $\mu', \mu'', \nu_3^3$, 
$\nu_3^4$, where:
\begin{gather*}
\mu'(\sigma_1) = s_1w_0\,, \quad \mu'(\sigma_2) = s_2w_0\,, \quad \mu'(\sigma_3) = s_3w_0\,;\\
\mu''(\sigma_1) = s_2s_1s_2 s_1s_3s_2 s_1s_2s_1 s_3s_2s_1 s_2s_3\,, \quad \mu''(\sigma_2) 
= s_3s_2s_1 s_2s_1s_3 s_2s_1s_2 s_3\,,\\ 
\mu''(\sigma_3) = s_1s_3\,;\\
\nu_3^3(\sigma_1) = s_2s_1s_2 s_1\,, \quad \nu_3^3(\sigma_2) = s_2s_1s_2 s_1s_3s_2 s_1s_2s_1 
s_3s_2s_1\,, \quad \nu_3^3(\sigma_3) = s_1s_2s_1 s_2\,;\\
\nu_3^4(\sigma_1) = s_3s_2s_1 s_2\,, \quad \nu_3^4(\sigma_2)= s_2s_3s_2 s_1\,, \quad 
\nu_3^4(\sigma_3) = s_2s_1s_2 s_3\,.
\end{gather*}}

\noindent
{\bf Remark.}
The up-to-center-epimorphisms $\mu'$ and $\mu''$ are ordinary while $\nu_3^3$ and $\nu_3^4$ are 
not.

\bigskip\noindent
{\bf Proof.}
The proof is a straightforward calculation that we have made with the package 
``Chevie'' of GAP. We first compute the set $X_1$ of triples $(x,y,z) \in W^3$ 
that satisfy
\[
xyxyx = yxyxy\,, \quad xz=zx\,, \quad yzy=zyz\,.
\]
We have $|X_1|=600$. Let $\sim$ be the equivalence relation on $X_1$ defined by $(x_1, y_1, z_1) 
\sim (x_2, y_2, z_2)$ if there exists $w \in W$ such that $x_2= wx_1w^{-1}$, $y_2=wy_1w^{-1}$, and 
$z_2=wz_1w^{-1}$. We compute a representative for each class, and we denote by $X_2$ the set of 
these representatives. We have $|X_2|=18$. Observe that there is a unique subgroup of $W$ of index 
2, which is the kernel of the epimorphism $\gamma: W \to C_2$ which sends $s_i$ to $-1$ for all 
$i=1,2,3$. So, a homomorphism $\f: A \to W$ is a proper up-to--center-epimorphism if and only if 
$|\Im \f|= {|W| \over 2} = 60$. For $(x,y,z) \in X_2$, we denote by $W(x,y,z)$ the subgroup of $W$ 
generated by $x,y,z$. Let $X_3$ be the set of triples $(x,y,z) \in X_2$ that satisfy 
$|W(x,y,z)|=60$. We compute $X_3$ and obtain
\begin{multline*}
X_3= \{ (\mu'(\sigma_1), \mu'(\sigma_2), \mu'(\sigma_3)),\ (\mu''(\sigma_1), \mu''(\sigma_2), 
\mu''(\sigma_3)),\\ (\nu_3^3(\sigma_1), \nu_3^3(\sigma_2), \nu_3^3(\sigma_3)),\ 
(\nu_3^4(\sigma_1), \nu_3^4(\sigma_2), \nu_3^4(\sigma_3))\}\,.\quad\square
\end{multline*}

\bigskip\noindent
{\bf Proposition 2.6.} {\it Assume $\Gamma=H_3$. Then there are precisely two equivalence classes 
of proper up-to-center-epimorphisms. These are represented by $\mu'$ and $\nu_3^3$.}

\bigskip\noindent
{\bf Proof.} Let $\alpha: W \to W$ be the automorphism defined by
\[
\alpha(s_1)= s_2s_1s_3 s_2s_3s_1 s_2s_1s_3 s_2s_3s_1 s_2\,, \quad \alpha(s_2)= s_2\,, \quad 
\alpha(s_3)=s_3\,.
\]
Let
\[
u= s_1s_2s_1 s_2s_1s_3 s_2\,, \quad v=s_1s_2s_1 s_3\,.
\]
Then $\mu''=\alpha \circ \Conj_u \circ \mu'$ and $\nu_3^4= \alpha \circ \Conj_v \circ \nu_3^3$. On 
the other hand, $\mu'$ and $\nu_3^3$ are not equivalent since $\mu'$ is ordinary while $\nu_3^3$ 
is not.
\endproof

\bigskip\noindent
{\bf Proposition 2.7.} {\it Assume $\Gamma=E_7$. Then there is a unique equivalence class of 
proper up-to-center-epimorphisms represented by $\mu': A \to W$, where
\[
\mu'(\sigma_i) = w_0s_i \quad \text{for } 1 \le i\le 7\,,
\]
and $w_0$ denotes the longest element of $W$. In particular, any proper up-to-center-epimorphism 
$A \to W$ is ordinary.}

\bigskip\noindent
{\bf Proof.} The proof is a straightforward but not easy calculation. It has been made with the 
package ``Chevie'' of GAP, using several computers at the same time.

\bigskip\noindent
First, observe that the homomorphism $\mu': A \to W$ defined in the statement of the lemma is a 
proper up-to-center-epimorphism. Thus, it suffices to show that there is a unique equivalence 
class of proper up-to-center-epimorphisms. Actually, we prove that there is a unique conjugacy class 
of proper up-to-center-epimorphisms, but this is not a contradiction as all the automorphisms of 
$W(E_7)$ are inner (see \cite{Fra}).

\bigskip\noindent
Let $\f: A \to W$ be a proper up-to-center-epimorphism, and write $t_i= \f(\sigma_i)$ for all $1 
\le i\le 7$. One can easily verify the following facts.
\begin{enumerate}
\item $t_i$ is conjugate to $t_1$ for all $1 \le i\le 7$;
\item $t_i \neq t_j$ for $i \neq j$;
\item Let $\gamma: W \to C_2$ be the epimorphism defined by $\gamma(s_i)=-1$ for all $1 \le i\le 
7$. Then $\Im \f = \Ker \gamma$.
\end{enumerate}
In particular, we have $\gamma(t_i)=1$ for all $1 \le i\le 7$.

\bigskip\noindent
The group $W=W(E_7)$ has 60 conjugacy classes that we number $C_1, C_2, \dots, C_{60}$. Call 
$C_j$ {\it even} if $\gamma(w)=1$ for all $w \in C_j$. Half of these classes, say $C_1, C_2, \dots, 
C_{30}$, are even. Moreover, we can assume that $C_1$ is the trivial class, $\{1\}$. So, we just 
need to consider the classes $C_2, C_3, \dots, C_{30}$.

\bigskip\noindent
For each $j \in \{2,3, \dots, 30\}$, we choose an element $T_{j\,4} \in C_j$ and we put
\begin{align*}
X_j=& \{u \in C_j\ ;\ u T_{j\,4} u= T_{j\,4} u T_{j\,4} \}\,,\\
Y_j=& \{ (u,v) \in (C_j)^2\ ;\ uv=vu\}\,,\\
Z_j=& \{ (u_2,u_3,u_5) \in (X_j)^3\ ;\ (u_r,u_s) \in Y_j \text{ for } r,s \in \{2,3,5\},r<s\}\,,\\
\tilde Z_j=& \{ (u_2,u_3,u_5) \in Z_j\ ;\ u_2 \neq u_3 \neq u_5 \neq u_2\}\,.
\end{align*}
For each $j \in \{2, \dots, 30\}$, we first compute $X_j$ and $Y_j$, then we compute $Z_j$ from 
$X_j$ and $Y_j$, and, finally, we compute $\tilde Z_j$ from $Z_j$. We give in Table 1 the 
cardinalities of $C_j, X_j, Y_j, Z_j$ and $\tilde Z_j$, for all $j \in \{2, \dots, 30\}$.
\[
\begin{array}{|c|c|c|c|c|c|c|c|c|c|c|c|}
\hline
j&2&3&4&5&6&7&8&9&10&11\\
\hline
C_j&63&315&945&3780&672&2240&13440&3780&7560&7560\\
\hline
X_j&33&129&161&225&136&136&208&241&129&129\\
\hline
Y_j&1953&13545&61425&196560&28224&58240&161280&120960&151200&151200\\
\hline
Z_j&4353&26625&14081&5505&3016&1000&352&9601&5217&5217\\
\hline
\tilde Z_j&2079&18432&8256&2112&720&216&36&5472&3120&3120\\
\hline
\hline
j&12&13&14&15&16&17&18&19&20&21\\
\hline
C_j&11340&45360&48384&10080&10080&20160&30240&40320&40320&120960\\
\hline
X_j&145&97&56&120&72&48&56&72&114&78\\
\hline
Y_j&362880&362880&167349&221760&141120&201600&302400&161280&80640&241920\\
\hline
Z_j&4993&1153&56&408&576&48&176&72&114&78\\
\hline
\tilde Z_j&2592&384&0&0&0&0&0&0&0&0\\
\hline
\hline
j&22&23&24&25&26&27&28&29&30&\\
\hline
C_j&207360&90720&90720&161280&145152&60480&60480&120960&96768&\\
\hline
X_j&57&57&89&37&58&48&48&64&41&\\
\hline
Y_j&1000397&725760&725760&797935&580608&219705&241920&483840&774144&\\
\hline
Z_j&57&513&545&37&58&108&108&64&41&\\
\hline
\tilde Z_j&0&144&144&0&0&0&0&0&0&\\
\hline
\end{array}\]

\centerline{{\bf Table 1.} Cardinalities of $C_j$, $X_j$, $Y_j$, $Z_j$, and $\tilde Z_j$.} 

\bigskip\noindent
If $\tilde Z_j= \emptyset$, then $t_1, \dots, t_7$ cannot lie in $C_j$. So, it suffices to 
consider those classes for which $\tilde Z_j \neq \emptyset$, that is, $j \in I=\{2,3, \dots, 13, 
23,24\}$ (see Table 1). For each $j \in I$ we compute the set
\[
U_j= \{ (u,v) \in (C_j)^2\ ;\ uvu=vuv\}\,,
\]
and, from $U_j$, $Y_j$ and $\tilde Z_j$, we compute the set $V_j$ of 7-tuples 
$(u_1,\dots,u_7) \in (C_j)^7$ that satisfy
\begin{gather*}
u_4=T_{j\,4}\,,\quad
(u_2,u_3,u_5) \in \tilde Z_j\,, \quad (u_1,u_i) \in Y_j \text{ for } i\in\{2,4,5,6,7\}\,,\\
(u_6,u_i) \in Y_j \text{ for }i\in \{2,3,4\}\,, \quad (u_7,u_i)\in Y_j \text{ for } i\in \{2,3,4,5\}\,,\\
(u_1,u_3), (u_5,u_6), (u_6,u_7) \in U_j\,.
\end{gather*} 
The cardinalities of $U_j$ and $V_j$ are given in Table 2. Observe that $V_j=\emptyset$ for 
all $j \in I$ except for $j=2$, thus we must have $t_1, t_2, \dots, t_7 \in C_2$.
\[
\begin{array}{|c|c|c|c|c|c|c|c|}
\hline
j&2&3&4&5&6&7&8\\
\hline
U_j&2079&40635&152145&850500&91392&304640&2795520\\
\hline
V_j&23040&0&0&0&0&0&0\\
\hline
\hline
j&9&10&11&12&13&23&24\\
\hline
U_j&910980&975240&975240&1644300&4399920&4000914&8074080\\
\hline
V_j&0&0&0&0&0&0&0\\
\hline
\end{array}
\]

\centerline{{\bf Table 2.} Cardinalities of $U_j$ and $V_j$.}

\bigskip\noindent
Consider the set
\[
V_2'= \{ (uw_0s_1u^{-1}, \dots, uw_0s_7u^{-1}); u \in W\text{ and } uw_0s_4u^{-1} = T_{2\,4}\}\,.
\]
We calculate $V_2'$ and observe that $|V_2'|=|V_2|=23040$. Thus, $V_2'=V_2$ and $\f$ must be 
conjugate to $\mu'$.
\endproof

\section{Composition}

The purpose of the present section is to show that an extraordinary up-to-center-epimorphism 
cannot be derived from an endomorphism of $A$. More precisely, we prove the following.

\bigskip\noindent
{\bf Theorem 3.1.} {\it Assume $\Gamma$ to be connected. Let $\Phi: A \to A$ be an endomorphism. 
If $\mu \circ \Phi: A \to W$ is an up-to-center-epimorphism, then $\mu \circ \Phi$ is ordinary. 
In particular, we have $\Phi(CA) \subset CA$ if $\mu \circ \Phi$ is an up-to-center-epimorphism.}

\bigskip\noindent
Obviously, there is nothing to prove if all the up-to-center-epimorphisms $A \to W$ are 
ordinary. So, by \cite{CoPa} and the previous section, we just need to study those Artin groups 
whose Coxeter graphs lie in the set $\{ I_2(p); p\ge 4 \text{ and } p\equiv 0(\mod\,4)\} \cup 
\{B_n; n\ge 3\} \cup \{D_n; n\ge 5 \text{ and } n\text{ odd}\} \cup \{A_3,H_3\}$.

\bigskip\noindent
We prove Theorem 3.1 with a case by case study, using the homomorphism $U: A \to \Sym( \Z \times 
T)$ introduced in \cite{CoPa} and defined as follows.

\bigskip\noindent
Let $T=\{ ws_iw^{-1}; 1\le i\le n\text{ and }w \in W\}$ be the {\it set of reflections} of $W$, 
and, for $1\le i\le n$, define the bijection $U_i: \Z \times T \to \Z \times T$ by
\[
U_i(k,t)= \left\{\begin{array}{ll}
(k,s_its_i)&\quad\text{if }t \neq s_i\,,\\
(k+1,t)&\quad\text{if }t=s_i\,.
\end{array}\right.
\]

\vfill\eject

\bigskip\noindent
{\bf Proposition 3.2 (Cohen, Paris \cite{CoPa}).}{\it
\begin{enumerate}
\item The map $\sigma_i \to U_i$, $1 \le i\le n$, determines a homomorphism $U: A \to
\Sym (\Z \times T)$.
\item Let $g \in A$ and $(k,t) \in \Z \times T$. Write $(k',t')= U(g)(k,t)$. Then $U(g)(k+l,t)= 
(k'+l,t')$ for all $l \in \Z$.
\item Let $g_1,g_2 \in G$ and $(k,t) \in \Z \times T$. Write $(k_1,t_1)=U(g_1)(k,t)$ and 
$(k_2,t_2)=U(g_2)(k,t)$. If $\mu(g_1)= \mu(g_2)$, then $t_1=t_2$ and $k_1 \equiv k_2 (\mod\,2)$.
\end{enumerate}}

\noindent
Now, we start the proof of Theorem 3.1 with the Artin groups of rank 2.

\bigskip\noindent
{\bf Lemma 3.3.} {\it Suppose $\Gamma=I_2(p)$ with $p \ge 4$ and $p \equiv 0 (\mod\,4)$. Let 
$\Phi: A \to A$ be an endomorphism. If $\mu \circ \Phi$ is an up-to-center-epimorphism, then $\mu 
\circ \Phi$ is ordinary.}

\bigskip\noindent
{\bf Proof.} Write $p=4q$, where $q \in \N$. Recall from Section 2 that any up-to-center-epimorphism 
$A \to W$ is an epimorphism. Recall also, from \cite{CoPa}, Proposition 2.2, 
that there are precisely two equivalence classes of extraordinary epimorphisms 
represented by $\nu_p^1, \nu_p^2: A \to W$, where
\[
\begin{array}{ll}
\nu_p^1(\sigma_1)=s_1\,, &\quad \nu_p^1(\sigma_2)=s_2s_1\,,\\
\nu_p^2(\sigma_1)=s_1s_2\,,&\quad \nu_p^2(\sigma_2)=s_2\,.
\end{array}
\]
Now, take and endomorphism $\Phi: A \to A$, and suppose that $\mu \circ \Phi$ is an (up-to-center) 
epimorphism.

\bigskip\noindent
{\bf Assertion 1.} {\it We have $\mu \circ \Phi \neq \nu_p^1$.}

\bigskip\noindent
{\bf Proof.} Suppose $\mu \circ \Phi = \nu_p^1$. Put
\[
\begin{array}{rcll}
t_{2i-1}&=&(s_2s_1)^{i-1} s_2 (s_1s_2)^{i-1} &\quad\text{for }1 \le i\le q\,,\\
t_{2i}&=&s_1(s_2s_1)^{i-1} s_2 (s_1s_2)^{i-1} s_1 &\quad\text{for }1\le i\le q\,.
\end{array}
\]
A direct calculation gives
\begin{align*}
U(\sigma_1):&\left\{
\begin{array}{llll}
(k,t_{2i-1})&\mapsto&(k,t_{2i})&\quad\text{for }1\le i\le q\,,\\
(k,t_{2i})&\mapsto&(k,t_{2i-1})&\quad\text{for }1\le i\le q\,.
\end{array}\right.\\
\\
U(\sigma_2\sigma_1):&\left\{
\begin{array}{llll}
(k,t_1)&\mapsto&(k,t_3)\,,\\
(k,t_2)&\mapsto&(k+1,t_1)\,,\\
(k,t_{2i-1})&\mapsto&(k,t_{2i+1})&\quad\text{for }2\le i\le q-1\,,\\
(k,t_{2i})&\mapsto&(k,t_{2i-2})&\quad\text{for }2\le i\le q-1\,,\\
(k,t_{2q-1})&\mapsto&(k,t_{2q})\,,\\
(k,t_{2q})&\mapsto&(k,t_{2q-2})\,.
\end{array}\right.
\end{align*}
Observe that $\mu(\sigma_1)=\nu_p^1(\sigma_1) =(\mu \circ \Phi) (\sigma_1)$ and $\mu(\sigma_2 
\sigma_1)= \nu_p^1(\sigma_2)= (\mu \circ \Phi)(\sigma_2)$, thus, by Proposition 3.2, there exist 
$a_1, \dots, a_{2q}, b_1, \dots, b_{2q} \in \Z$ such that
\begin{align*}
U(\Phi(\sigma_1)):&\left\{
\begin{array}{llll}
(k,t_{2i-1})&\mapsto&(k+2a_{2i-1},t_{2i})&\quad\text{for }1\le i\le q\,,\\
(k,t_{2i})&\mapsto&(k+2a_{2i},t_{2i-1})&\quad\text{for }1\le i\le q\,.
\end{array}\right.\\
\\
U(\Phi(\sigma_2)):&\left\{
\begin{array}{llll}
(k,t_1)&\mapsto&(k+2b_1,t_3)\,,\\
(k,t_2)&\mapsto&(k+2b_2+1,t_1)\,,\\
(k,t_{2i-1})&\mapsto&(k+2b_{2i-1},t_{2i+1})&\quad\text{for }2\le i\le q-1\,,\\
(k,t_{2i})&\mapsto&(k+2b_{2i},t_{2i-2})&\quad\text{for }2\le i\le q-1\,,\\
(k,t_{2q-1})&\mapsto&(k+2b_{2q-1},t_{2q})\,,\\
(k,t_{2q})&\mapsto&(k+2b_{2q},t_{2q-2})\,.
\end{array}\right.
\end{align*}
With a direct calculation we obtain
\[
U(\Phi(\sigma_1\sigma_2)^{2q}):\left\{
\begin{array}{lll}
(0,t_1)&\mapsto&(2q(a_2+a_3+b_1+b_4),t_1)\,,\\
(0,t_2)&\mapsto&(4q(a_1+b_2)+2q,t_2)\,,\\
(0,t_{2i-1})&\mapsto&(2q(a_{2i}+a_{2i+1}+b_{2i-1}+b_{2i+2}),t_{2i-1})\\
&&\quad\text{for }2\le i\le q-1\,,\\
(0,t_{2i})&\mapsto&(2q(a_{2i-2}+a_{2i-1}+b_{2i-3}+b_{2i}),t_{2i})\\
&&\quad\text{for }2\le i\le q-1\,,\\
(0,t_{2q-1})&\mapsto&(4q(a_{2q}+b_{2q-1}),t_{2q-1})\,,\\
(0,t_{2q})&\mapsto&(2q(a_{2q-2}+a_{2q-1}+b_{2q-3}+b_{2q}),t_{2q})\,.
\end{array}\right.\\
\]\[
U(\Phi(\sigma_2\sigma_1)^{2q}):\left\{
\begin{array}{llll}
(0,t_1)&\mapsto&(4q(a_1+b_2)+2q,t_1)\,,\\
(0,t_{2i})&\mapsto&(2q(a_{2i}+a_{2i+1}+b_{2i-1}+b_{2i+2}),t_{2i})\\
&&\quad\text{for }1\le i\le q-1\,,\\
(0,t_{2i-1})&\mapsto&(2q(a_{2i-2}+a_{2i-1}+b_{2i-3}+b_{2i}),t_{2i-1})\\
&&\quad\text{for }2\le i\le q\,,\\
(0,t_{2q})&\mapsto&(4q(a_{2q}+b_{2q-1}),t_{2q})\,.
\end{array}\right.
\]
The equalities $U(\Phi(\sigma_1\sigma_2)^{2q})(0,t_j) = U(\Phi(\sigma_2\sigma_1)^{2q})(0,t_j)$, for 
$j=1, \dots, 2q$, are equivalent 
to the equations
\begin{gather*}
(a_2+a_3+b_1+b_4)-2(a_1+b_2)=1\,,\\
(a_{2i}+a_{2i+1}+b_{2i-1}+b_{2i+2})-(a_{2i-2}+a_{2i-1}+b_{2i-3}+b_{2i})=0\quad\text{for }2\le i\le q-1\,,\\
2(a_{2q}+b_{2q-1})-(a_{2q-2}+a_{2q-1}+b_{2q-3}+b_{2q})=0\,.
\end{gather*}
The sum of these equations gives $2(a_{2q}+b_{2q-1}-a_1-b_2)=1$, which is clearly impossible.

\bigskip\noindent
{\bf Assertion 2.} {\it Let $\alpha: W \to W$ be an automorphism. Then there exists an 
endomorphism $\tilde \alpha: A \to A$ (not necessarily bijective) such that $\alpha \circ \mu = 
\mu \circ \tilde \alpha$.}

\bigskip\noindent
{\bf Proof.} Let $\alpha_0: W \to W$ denote the automorphism determined by $\alpha_0(s_1)=s_2$ 
and $\alpha_0(s_2)=s_1$. For $k \in \{ 0,1, \dots, 2q-1\}$ with $\gcd(2k+1,2q)=1$, let 
$\alpha_k: W \to W$ denote the automorphism determined by $\alpha_k(s_1)=s_1$ and $\alpha_k(s_2)= 
(s_2s_1)^k s_2 (s_1s_2)^k$. It is easily checked that $\Aut (W)$ is generated as a monoid by the 
set $\{\Conj_w; w\in W\} \cup \{\alpha_0\} \cup \{\alpha_k; 0\le k\le q-1\text{ and 
}\gcd(2k+1,4q)=1\}$. So, it suffices to prove Assertion 2 for $\alpha$ in this set of generators.

\bigskip\noindent
Let $w \in W$. Take $u \in A$ such that $\mu(u)=w$. Then $\Conj_w \circ \mu= \mu \circ 
\Conj_u$. Let $\tilde \alpha_0: A \to A$ be the automorphism determined by $\tilde 
\alpha_0(\sigma_1) = \sigma_2$ and $\tilde \alpha_0 (\sigma_2)= \sigma_1$. Then $\alpha_0 \circ 
\mu= \mu \circ \tilde \alpha_0$. Let $k \in \{0,1, \dots, q-1\}$ with $\gcd(2k+1,2q)=1$, and let 
$\tilde \alpha_k: A \to A$ be the endomorphism defined by $\tilde \alpha_k (\sigma_1) = \sigma_1$ 
and $\tilde \alpha_k(\sigma_2)= (\sigma_2\sigma_1)^k \sigma_2 (\sigma_1 \sigma_2)^k$. Then 
$\alpha_k \circ \mu= \mu \circ \tilde \alpha_k$.

\bigskip\noindent
{\bf End of the proof of Lemma 3.3.} Suppose that $\mu \circ \Phi$ is extraordinary. By 
\cite{CoPa}, Proposition 2.2, $\mu \circ \Phi$ must be equivalent to either $\nu_p^1$ or 
$\nu_p^2$. Without loss of generality, we may assume that $\mu \circ \Phi$ is equivalent to 
$\nu_p^1$. Take an automorphism $\alpha: W \to W$ such that $\alpha \circ \mu \circ \Phi= 
\nu_p^1$. By Assertion 2, one can find an endomorphism $\tilde \alpha: A \to A$ such that $\alpha 
\circ \mu= \mu \circ \tilde \alpha$. Then we have the equality $\mu \circ (\tilde \alpha \circ 
\Phi)= \nu_p^1$ which contradicts Assertion 1.
\endproof

\bigskip\noindent
{\bf Lemma 3.4.} {\it Suppose $\Gamma=B_n$ with $n \ge 3$. Let $\Phi: A \to A$ be an 
endomorphism. If $\mu \circ \Phi$ is an up-to-center-epimorphism, then $\mu \circ \Phi$ is 
ordinary.}

\bigskip\noindent
{\bf Proof.} By \cite{CoPa}, Proposition 3.2, and by Proposition 2.4 of the present paper, we have 
the following.
\begin{itemize}
\item If $n$ is even, then there is precisely one equivalence class of extraordinary 
up-to-center-epimorphisms represented by the epimorphism $\nu_n^1$ of \cite{CoPa}, Proposition 3.2.
\item If $n$ is odd, then there are precisely four equivalence classes of extraordinary 
up-to-center-epimorphisms. These are 
represented by the epimorphisms $\nu_n^1, \nu_n^2$ of \cite{CoPa}, 
Proposition 3.2, and by the proper up-to-center-epimorphisms $\nu_n^3,\nu_n^4$ of Proposition 2.4.
\end{itemize}
Recall that, for $j=1, \dots, n$, $e_j$ denotes the element $s_j \dots s_2s_1s_2 \dots s_j$ of 
$W$. Observe that, for $i=1,2,3,4$, we have $\nu_n^i(\sigma_j)= e_{j-1} s_j$ for all $2 \le j\le 
n$. A careful reading of the proof of \cite{CoPa}, Lemma 3.10, shows that, 
for a given endomorphism $\Phi: A \to A$, we cannot have $\mu \circ \Phi(\sigma_2)= e_1s_2$ and 
$\mu \circ \Phi (\sigma_3)= e_2s_3$. Thus $\mu \circ \Phi \neq \nu_n^i$ for all 
$i=1,2,3,4$.

\bigskip\noindent
It is shown in the proof of \cite{CoPa}, Proposition 3.9, that, for any automorphism $\alpha:W \to W$, 
there exists an endomorphism $\tilde \alpha: A \to A$ such that $\alpha \circ \mu= \mu \circ 
\tilde \alpha$. Now, suppose that $\mu \circ \Phi$ is extraordinary. Then there exists $i \in 
\{1,2,3,4\}$ and $\alpha \in \Aut(W)$ such that $\alpha \circ \mu \circ \Phi= \nu_n^i$. Take 
$\tilde \alpha:A \to A$ such that $\alpha \circ \mu= \mu \circ \tilde \alpha$. Then $\mu \circ 
(\tilde \alpha \circ \Phi)= \nu_n^i$: a contradiction.
\endproof

\bigskip\noindent
{\bf Lemma 3.5.} {\it Suppose $\Gamma=A_3$, or $\Gamma=D_n$, with $n\ge 5$ and $n$ odd. Let 
$\Phi: A \to A$ be an endomorphism. If $\mu \circ \Phi$ is an up-to-center-epimorphism, then $\mu 
\circ \Phi$ is ordinary.}

\bigskip\noindent
{\bf Proof.} We have shown in Section 2 that there is no proper up-to-center-epimorphism $A \to 
W$. Furthermore, the hypothesis ``$\Phi$ is an epimorphism'' is not needed in the proof of 
\cite{CoPa}, Lemma 4.6. Therefore, Lemma 3.5 follows from \cite{CoPa}, Lemma 4.6, together with 
the fact that all the automorphisms of $W$ are inner (see \cite{Fra}).
\endproof

\bigskip\noindent
Now, the next lemma finishes the proof of Theorem 3.1.

\bigskip\noindent
{\bf Lemma 3.6.} {\it Assume $\Gamma=H_3$. Let $\Phi: A \to A$ be an endomorphism. If $\mu \circ 
\Phi$ is an up-to-center-epimorphism, then $\mu \circ \Phi$ is ordinary.}

\bigskip\noindent
{\bf Proof.} We follow the same strategy as in the proof of \cite{CoPa}, Proposition 5.3. By 
\cite{CoPa}, Lemma 5.1, and by Lemma 2.5 in the present paper, there are precisely four conjugacy 
classes of up-to-center-epimorphisms. These are represented by the epimorphisms $\nu_3^1, 
\nu_3^2$ of \cite{CoPa}, Lemma 5.1, and by the proper up-to-center-epimorphisms $\nu_3^3, 
\nu_3^4$ of Lemma 2.5.

\bigskip\noindent
{\bf Assertion 1.} {\it Let $\Phi: A \to A$ be an endomorphism. Then $\mu \circ \Phi \neq 
\nu_3^i$ for all $i=1,2,3,4$.}

\bigskip\noindent
{\bf Proof.} Let
\begin{gather*}
r_1=s_1\,, \quad r_2=s_2\,, \quad r_3=s_3\,, \quad r_4=s_1s_2s_1\,, \quad r_5=s_2s_3s_2\,, \quad 
r_6=s_2s_1s_2\,, \\
r_7=s_1s_2s_3 s_2s_1\,, \quad r_8=s_3s_2s_1 s_2s_3\,, \quad r_9=s_1s_2s_1 s_2s_1\,,\quad 
r_{10}=s_1s_3s_2 s_1s_2s_3 s_1\,, \\
\quad r_{11}=s_2s_1s_2 s_3s_2s_1 s_2\,, \quad r_{12}=s_2s_1s_3 s_2s_1s_2 s_3s_1s_2\,,\quad
r_{13}=s_1s_2s_1 s_2s_3s_2 s_1s_2s_1\,, \\
\quad r_{14}=s_1s_2s_1 s_3s_2s_1 s_2s_3s_1 s_2s_1\,, 
\quad r_{15}=s_2s_1s_2 s_1s_3s_2 s_1s_2s_3 s_1s_2s_1 s_2\,.
\end{gather*}
Suppose $\mu \circ \Phi = \nu_3^1$. It is shown in the proof of \cite{CoPa}, Proposition 5.3, 
that the equalities $U( \Phi( \sigma_2 \sigma_3 \sigma_2))(0,r_j) = U(\Phi (\sigma_3 \sigma_2 
\sigma_3)) (0,r_j)$ for $j=2,5,6,8,10,11$ are equivalent to the equations
\begin{align*}
b_2-b_{10}+b_{12}-c_2-c_3+c_5&=-1\,,\\
b_5-b_{12}+b_{14}-c_2-c_5+c_{11}&=0\,,\\
b_5+b_6-b_7+c_4-c_6-c_8&=-1\,,\\
b_7+b_8-b_{11}-c_1+c_6-c_8&=0\,,\\
b_2-b_6+b_{10}+c_3-c_4-c_{10}&=0\,,\\
b_8+b_{11}-b_{14}+c_1-c_{10}-c_{11}&=-1\,,
\end{align*}
where $b_1,b_2, \dots, b_{15}, c_1, c_2, \dots, c_{15} \in \Z$. The sum of these 6 equations 
gives $2(b_2+b_5+b_8-c_2-c_8-c_{10})=-3$ which clearly cannot hold.

\bigskip\noindent
It is also shown in the proof of \cite{CoPa}, Proposition 5.3, that $\mu \circ \Phi \neq 
\nu_3^2$. So, it remains to show that $\mu \circ \Phi \neq \nu_3^3$ and $\mu \circ \Phi \neq 
\nu_3^4$. We proceed with the same method.

\bigskip\noindent
Suppose $\mu \circ \Phi = \nu_3^3$. Then the equalities $U(\Phi (\sigma_2 \sigma_3 \sigma_2)) 
(0,r_j) = U(\Phi (\sigma_3 \sigma_2 \sigma_3)) (0,r_j)$ for $j=2,4,6, 9,12,13$ are equivalent to 
the equations
\begin{align*}
b_2-b_4+b_{11}-c_1-c_2+c_9&=0\,,\\
b_2+b_4-b_7+c_1-c_4-c_6&=1\,,\\
b_6-b_9+b_{13}+c_3-c_6-c_{14}&=0\,,\\
b_9-b_{11}+b_{12}-c_2-c_9+c_{14}&=0\,,\\
b_7-b_{10}+b_{12}+c_4-c_{12}-c_{13}&=-1\,,\\
-b_6+b_{10}+b_{13}-c_3+c_{12}-c_{13}&=-1\,,
\end{align*}
where $b_1,b_2, \dots, b_{15}, c_1, c_2, \dots, c_{15} \in \Z$. The sum of these 6 equations 
gives $2(b_2+b_{12}+b_{13}-c_2-c_6-c_{13})=-1$ which clearly cannot hold.

\bigskip\noindent
Suppose $\mu \circ \Phi= \nu_3^4$. Then the equalities $U( \Phi (\sigma_2 \sigma_3 \sigma_2)) 
(0,r_j)= U( \Phi( \sigma_3 \sigma_2 \sigma_3)) (0,r_j)$ for $j=4,7,8,9,13,14$ are equivalent to 
the equations
\begin{align*}
b_4+b_8-b_{12}-c_4+c_{11}-c_{14}&=0\,,\\
b_7+b_{12}-b_{15}+c_4-c_7-c_{13}&=0\,,\\
b_1-b_4+b_8-c_8+c_9-c_{11}&=0\,,\\
-b_1+b_7+b_9+c_6-c_8-c_9&=0\,,\\
b_{13}-b_{14}+b_{15}-c_5+c_7-c_{13}&=0\,,\\
-b_9+b_{13}+b_{14}+c_5-c_6-c_{14}&=1\,,
\end{align*}
where $b_1,b_2, \dots, b_{15},c_1,c_2, \dots, c_{15} \in \Z$. The sum of these equations gives 
$2(b_7+b_8+b_{13}-c_8-c_{13}-c_{14})=1$ which clearly cannot hold.

\bigskip\noindent
{\bf End of the proof of Lemma 3.6.} Let $\Phi: A \to A$ be an endomorphism, and suppose that 
$\mu \circ \Phi$ is an extraordinary up-to-center-epimorphism. Then there exist $i \in \{1,2,3,4 
\}$ and $w \in W$ such that $\Conj_w \circ \mu \circ \Phi = \nu_3^i$. Take $u \in A$ such that 
$\mu(u)=w$. Then we have the equality $\mu \circ (\Conj_u \circ \Phi)= \nu_3^i$ which cannot hold 
by Assertion 1.
\endproof


\section{Coloured Artin groups}

We turn now to prove our main result.

\bigskip\noindent
{\bf Theorem 4.1.} {\it Let $A$ be a (non necessarily irreducible) spherical type Artin group. 
Then the coloured Artin group $CA$ is a characteristic subgroup of $A$.}

\bigskip\noindent
The following proposition is a preliminary result to the proof of Theorem 4.1. Part 1 can be found 
in \cite{Par}, Proposition 4.2, and Part 2 is a mixture of \cite{Par}, Proposition 4.3 and 
\cite{Par}, Proposition 5.1.

\bigskip\noindent
First, recall that, if $H_1,H_2$ are two subgroups of a given group $G$, then $(H_1, H_2)$ 
denotes the subgroup of $G$ generated by $\{x_1^{-1} x_2^{-1} x_1 x_2; x_1 \in H_1 \text{ and } 
x_2 \in H_2\}$. In particular, the equality $(H_1,H_2)=\{1\}$ means that every element of $H_1$ 
commutes with every element of $H_2$.

\bigskip\noindent
{\bf Proposition 4.2 (Paris \cite{Par}).} {\it
\begin{enumerate}
\item Assume $\Gamma$ to be connected. Let $H_1, H_2$ be two subgroups of $A$ such that 
$(H_1,H_2)=\{1\}$ and $A=H_1 \cdot H_2$. Then either $H_1 \subset Z(A)$ or $H_2 \subset Z(A)$.
\item Let $\Gamma, \Omega$ be two connected spherical type Coxeter graphs. Suppose that there 
exists a monomorphism $\Phi: A (\Gamma) \to A(\Omega)$ such that $A(\Omega)= \Im \Phi \cdot Z( A( 
\Omega))$. Then $\Gamma \simeq \Omega$.
\end{enumerate}}

\noindent
{\bf Proof of Theorem 4.1.} Let $\Gamma$ be the Coxeter graph of $A$. Let $\Gamma_1, \dots, 
\Gamma_p$ be the connected components of $\Gamma$. For $j=1, \dots, p$, we denote by $\mu_j: 
A(\Gamma_j) \to W(\Gamma_j)$ the standard epimorphism. Note that $A= A(\Gamma_1) \times \dots 
\times A(\Gamma_p)$, $W= W(\Gamma_1) \times \dots \times W(\Gamma_p)$, and the standard 
epimorphism $\mu: A \to W$ is given by $\mu(g_1, \dots ,g_p)= (\mu_1(g_1), \dots, 
\mu_p(g_p))$. Up to a permutation of the $\Gamma_j$'s, we can assume that there is an $x \in 
\{0,1, \dots, p\}$ such that $\Gamma_j \neq A_1$ for all $1 \le j\le x$, and $\Gamma_j=A_1$ for 
all $x+1\le j\le p$. In particular, $A(\Gamma_1), \dots, A(\Gamma_x)$ are non-abelian irreducible 
Artin groups, and $A(\Gamma_{x+1}), \dots, A(\Gamma_p)$ are all isomorphic to $\Z$.

\bigskip\noindent
Let fix an isomorphism $\Phi: A \to A$. For $1 \le j\le p$, we denote by $\iota_j: A(\Gamma_j) 
\to A$ the natural embedding, and by $\kappa_j: A \to A(\Gamma_j)$ the projection on the $j$-th 
component. For $1 \le j,k \le p$, we put $\Phi_{j\,k}= \kappa_k \circ \Phi \circ \iota_j: 
A(\Gamma_j) \to A(\Gamma_k)$.

\bigskip\noindent
Let $k \in \{1, \dots, x\}$. Observe that $A(\Gamma_k) = \prod_{j=1}^p \Phi_{j\,k} (A(\Gamma_j))$ 
and $(\Phi_{j\,k} (A(\Gamma_j)), \Phi_{l\,k} (A(\Gamma_l)))$ $=\{1\}$ for $j \neq l$. Thus, by 
Proposition 4.2, there exists $\xi(k) \in \{1,2, \dots, p\}$ such that:
\begin{itemize}
\item $\Phi_{ \xi(k)\, k} (A( \Gamma_{\xi(k)})) \cdot Z( A(\Gamma_k)) = A(\Gamma_k)$;
\item $\Phi_{j\,k}( A(\Gamma_j)) \subset Z(A(\Gamma_k))$ if $j \neq \xi(k)$.
\end{itemize}
Furthermore, since $A(\Gamma_k)$ is non-abelian, the number $\xi(k)$ is unique and $\xi(k)\in 
\{1,2, \dots,$ $ x\}$.

\bigskip\noindent
Now, we prove that $\xi: \{1, \dots, x\} \to \{1, \dots, x\}$ is a surjection (and, therefore, a 
permutation). Suppose not, that is, there exists some $j \in \{1, \dots, x\}$ such that $\xi(k) 
\neq j$ for all $k \in \{1, \dots, x\}$. Then $\Phi_{j\,k}( A(\Gamma_j)) \subset Z(A(\Gamma_k))$ 
for all $k \in \{1, \dots, x\}$, thus $\Phi(A(\Gamma_j)) \subset Z(A)$. This contradicts the fact 
that $\Phi$ is an isomorphism and that $A(\Gamma_j)$ is non-abelian.

\bigskip\noindent
Now, we prove that $\Phi_{ \xi(k)\, k}: A(\Gamma_{\xi(k)}) \to A(\Gamma_k)$ is injective. Let $g 
\in \Ker \Phi_{ \xi(k)\,k}$. We have $\Phi_{\xi(k)\,l} (g) \in Z(A(\Gamma_l))$ for $l \neq k$ 
(since $\xi$ is a permutation), and $\Phi_{\xi(k)\,k}(g)=1$, thus $\Phi(g) \in Z(A)$. Since $\Phi$ 
is an isomorphism, it follows that $g \in Z(A) \cap A(\Gamma_{\xi(k)})= Z(A(\Gamma_{\xi(k)}))$. 
Recall from Proposition 2.1 that $Z(A(\Gamma_{\xi(k)}))$ is the infinite cyclic subgroup of 
$A(\Gamma_{\xi(k)})$ generated by the so-called standard generator, $\delta_{\xi(k)}$. So, there 
exists $l \in \Z$ such that $g= (\delta_{\xi(k)})^l$. Let $i_1, i_2, \dots, i_r$ be the vertices of 
$\Gamma_{\xi(k)}$ (which can be chosen in any order), and let $\pi= \sigma_{i_1} \sigma_{i_2} 
\dots \sigma_{i_r}$. It is shown in \cite{BrSa} that $\pi$ is not central in $A(\Gamma_{\xi(k)})$ 
(thus, by the above, $\pi \not\in \Ker \Phi_{\xi(k)\,k}$), but there exists some number $h>0$ 
such that $\delta_{\xi(k)}= \pi^h$. (This equality holds for any choice of the linear ordering of 
the vertices of $\Gamma_{\xi(k)}$.) The group $A(\Gamma_{\xi(k)})$ is torsion-free (see 
\cite{Del}), $\Phi_{\xi(k)\,k} (\pi) \neq 1$, and $\Phi_{\xi(k)\,k}(\pi)^{hl} = \Phi_{\xi(k)\,k} 
(g)=1$, thus $hl=0$, therefore $g=\pi^{hl}=1$.

\bigskip\noindent
In order to prove Theorem 4.1, we must show that $\Phi(CA) \subset CA$, that is, $(\mu \circ 
\Phi) (\sigma_i^2)=1$ for all $i=1, \dots, n$. Note that the inclusion $\Phi(CA) \subset CA$ 
implies that $\Phi(CA)=CA$ since $CA$ is of finite index in $A$ (the quotient $A/CA$ is isomorphic to $W$ 
which is finite). Take $i\in \{1, \dots, n\}$, and let $j \in \{1, \dots, p\}$ such that $i$ is a 
vertex of $\Gamma_j$. In order to prove that $(\mu \circ \Phi) (\sigma_i^2)=1$, it clearly 
suffices to show that $(\mu_k \circ \Phi_{j\,k}) (\sigma_i^2)=1$ for all $k \in \{1, \dots, p\}$.

\bigskip\noindent
Assume $j \in \{x+1, \dots, p\}$. In particular, $\Gamma_j \simeq A_1$ and $\Z \simeq A(\Gamma_j) 
\subset Z(A)$. Let $k \in \{1, \dots, p\}$. We have $\sigma_i \in Z(A)$, thus $\Phi(\sigma_i) \in 
Z(A)$, therefore $\Phi_{j\,k} (\sigma_i) \in Z(A(\Gamma_k))$. Let $\delta_k$ denote the standard 
generator of $Z(A(\Gamma_k))$. Then there exists $l \in \Z$ such that $\Phi_{j\,k}(\sigma_i)= 
\delta_k^l$. It is easily checked from the definition of $\delta_k$ (see Section 2) that 
$\mu_k(\delta_k^2)=1$. Hence, $(\mu_k \circ \Phi_{j\,k}) (\sigma_i^2)= \mu_k( \delta_k^{2l}) =1$.

\bigskip\noindent
Assume $j \in \{1, \dots, x\}$. Let $k \in \{1, \dots, p\}$. If $j \neq \xi(k)$, then 
$\Phi_{j\,k}(\sigma_i) \in Z(A(\Gamma_k))$, therefore, as in the previous case, $(\mu_k \circ 
\Phi_{j\,k})(\sigma_i^2)=1$. So, we can assume that $j = \xi(k)$. We have $\Im \Phi_{j\,k} \cdot 
Z(A(\Gamma_k)) = A(\Gamma_k)$ and $\Phi_{j\,k}: A(\Gamma_j) \to A(\Gamma_k)$ is injective, thus, 
by Proposition 4.2, $\Gamma_j \simeq \Gamma_k$. Moreover, the equality $\Im \Phi_{j\,k} \cdot 
Z(A(\Gamma_k)) = A(\Gamma_k)$ implies that $\Im (\mu_k \circ \Phi_{j\,k}) \cdot Z(W(\Gamma_k)) = 
W(\Gamma_k)$, namely, that $(\mu_k \circ \Phi_{j\,k}): A(\Gamma_j) \to W(\Gamma_k)$ is an 
up-to-center-epimorphism. We conclude by Theorem 3.1 that $\mu_k \circ \Phi_{j\,k}$ is ordinary and, 
consequently, that $(\mu_k \circ \Phi_{j\,k})(\sigma_i^2)=1$.
\endproof



\bigskip\bigskip\noindent
{\bf Nuno Franco},

\smallskip\noindent
R. Rom\~ao Ramalho 59,
Departamento de Matem\'atica, CIMA-UE,
Universidade de \'Evora,
7000 \'Evora,
Portugal

\smallskip\noindent
Institut de Math\'ematiques de Bourgogne, UMR 5584 du CNRS, Universit\'e de Bourgogne, B.P. 
47870, 21078 Dijon cedex, France

\smallskip\noindent
E-mail: {\tt nmf@evunix.uevora.pt} 

\bigskip\noindent
{\bf Luis Paris},

\smallskip\noindent 
Institut de Math\'ematiques de Bourgogne, UMR 5584 du CNRS, Universit\'e de Bourgogne, B.P. 
47870, 21078 Dijon cedex, France

\smallskip\noindent
E-mail: {\tt lparis@u-bourgogne.fr}

\end{document}